\documentclass{amsart}
\usepackage{amssymb,GriJer}

\usepackage{graphics}

\begin{document}
\title{Asymptotics of eigenfunctions on plane domains}

\author{Daniel Grieser}
\address{Institut f\"ur Mathematik, Carl-von-Ossietzky Universit\"at Oldenburg}
\email{grieser@mathematik.uni-oldenburg.de}
\author{David Jerison}
\address{Department of Mathematics, Massachusetts Institute of Technology}
\email{jerison@math.mit.edu}
\begin{abstract}
We consider a family of domains $(\Omega_N)_{N>0}$ obtained by attaching an $N\times 1$ rectangle
to a fixed set $\Omega_0 = \{(x,y):\, 0<y<1,\,-\phi(y)<x<0\}$, for a Lipschitz function $\phi\geq 0$.
We derive full asymptotic expansions, as $N\to\infty$, for the $m$th Dirichlet eigenvalue (for any fixed $m\in\NN$) and for the associated eigenfunction on $\Omega_N$. The second term involves a scattering phase arising in the Dirichlet problem on the infinite domain $\Omega_\infty$. We determine the first variation of this scattering phase, with respect to $\phi$,  at $\phi\equiv 0$. This is then used to prove sharpness of results, obtained previously by the same authors, about the location of extrema and nodal line of eigenfunctions on convex domains.
\end{abstract}

\maketitle


\tableofcontents

\section{Introduction}
\begin{figure}
\includegraphics{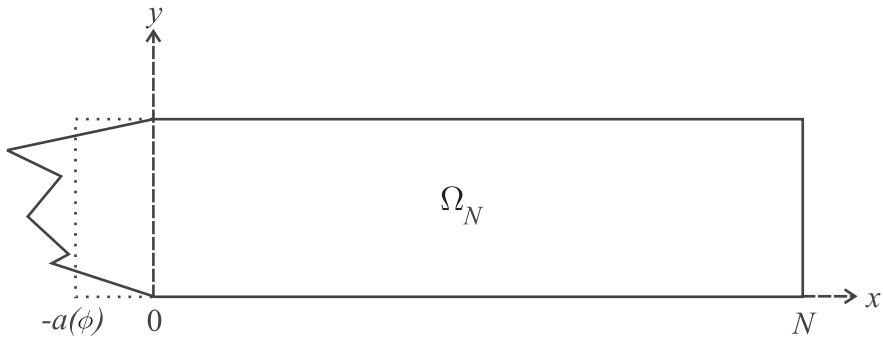}
\caption{The domain $\Omega_N$}
\label{figure1}
\end{figure}

For a Lipschitz function $\phi:[0,1]\to[0,\infty)$ and for $N\in[0,\infty]$ consider the plane domain (see Figure \ref{figure1})
\begin{equation}\label{GJ.0}
\Omega_N = \{(x,y)\in\bbR^2: 0<y<1,\, -\phi(y)<x<N\}
\end{equation}
and the eigenvalue problem for the Diríchlet Laplacian on $\Omega_N$:
\begin{align*}
(\Delta+\mu)u &= 0 \quad \text{on } \Omega_N\\
u &=0 \quad\text{at } \partial \Omega_N.
\end{align*}
For $N<\infty$ let $\mu_1(\Omega_N)<\mu_2(\Omega_N)\leq\dots$ be the eigenvalues, counted with multiplicities. The object of this paper is to study, for fixed $m\in\NN$, the asymptotic behavior of $\mu_m(\Omega_N)$ and of the associated eigenfunctions, as $N\to\infty$. This will then be used to answer some questions left open in our study (see \cite{Jer:DFNLCD,GriJer:AFNLCD,GriJer:SFECPD}) of the first and the second eigenfunction on general plane convex domains.

Our first main theorem is:

\begin{theorem}\label{GJ.4}
There is a number $a(\phi)\in[0,\max\phi]$ such that for each $m\in\NN$
the $m$th Dirichlet eigenvalue of $\Omega_N$ satisfies
\begin{equation}\label{GJ.5}
\mu_m(\Omega_N) = \pi^2 + \frac{m^2\pi^2}{(N+a(\phi))^2}
+ O(N^{-5}),\quad N\to\infty.
\end{equation}
In particular, the eigenvalues $\mu_1,\dots,\mu_m$ of $\Omega_N$ are simple for $N$
sufficiently large.
The suitably rescaled eigenfunction $u_m$ satisfies, for all multiindices
$\alpha$,
\begin{equation}\label{GJ.6}
\sup_{x>3\log N} |D^\alpha \left(u_m(x,y)-\sin m\pi\frac{x+a(\phi)}{N+a(\phi)}\,\sin\pi y\right)| = O(N^{-3}),
\end{equation}
\begin{equation}\label{GJ.6a}
\sup_{x\leq 3\log N} |u_m(x,y)| = O(N^{-1}{\log N}).
\end{equation}
The constants in the error terms only depend on $k,\ \alpha$ and  $\max\phi$.
\end{theorem}
Thus, the spectral data are very close to the data obtained on the rectangle $[-a(\phi),N]\times[0,1]$.
In fact, we will get complete asymptotic expansions for the eigenvalue and the eigenfunction and also much more precise information about the eigenfunction for small values of $x$, see Remark \ref{remnextterm} and Theorem \ref{ansatzcorrect} (in connection with \eqref{glue}).

The number $a(\phi)$ is closely related to a scattering phase in an associated non-compact problem, see Section \ref{secinf}, in particular Remark \ref{remarkscattering}. Therefore, its dependence on $\phi$ is very subtle.
Our second main theorem gives a perturbation analysis of $a(\phi)$ around
$\phi\equiv 0$:

\begin{theorem}\label{GJ.10}
Fix a Lipschitz function $\phi:[0,1]\to[0,\infty)$. Then, as $\epsilon \to 0$,
\begin{equation}\label{GJ.11}
a(\epsilon\phi) = 2\epsilon \int_0^1 \phi(y)(\sin \pi y)^2\,dy + O(\epsilon^2).
\end{equation}
\end{theorem}

As already mentioned, one motivation for the present study is to complement the results
about the first and second eigenfunction $u_1$, $u_2$ on a plane convex domain $\Omega$ obtained in \cite{Jer:DFNLCD,GriJer:SFECPD}. In these papers we considered the maximum of $u_1$ (which is w.l.o.g.\ assumed positive) and the nodal line of $u_2$, and a central goal was to localize these objects in terms of corresponding objects for eigenfunctions of an associated {\em ordinary} differential
operator. To state this more precisely, we first normalize $\Omega$ by a rotation and a dilation so that among all projections
of $\Omega$ onto lines the projection onto the $y$-axis has shortest possible length and this length equals one.
Let $\pi(x,y)=x$ be the projection map to the $x$ axis.
Let $I=\pi(\Omega)$  and $h:I\to(0,\infty)$ be the 'height function' of $\Omega$, that is,
$$h(x) = \text{length of the interval }\pi^{-1}(\{x\})\cap\Omega.$$
Let $\phi_1,\phi_2$ be the first and second eigenfunction of the Schr\"odinger operator
$-\dfrac{d^2}{dx^2} + \dfrac{\pi^2}{h(x)^2}$ on $I$, with Dirichlet boundary conditions (defined in terms of the variational principle on $H^1_0(I)$).
\begin{theorem}[\cite{Jer:DFNLCD,GriJer:SFECPD}]\label{GJ.20a}
Let the domain $\Omega$ be normalized as above and let $h:I\to(0,\infty)$ be its height function.
Denote by $\cM$ the set where $u_1$ achieves its maximum, by $\cN=u_2^{-1}(0)$ the nodal line, and let
$\{m\}$ and $\{n\}$ be the corresponding sets for $\phi_1,\phi_2$.

There is an absolute constant $C$ so that
\begin{equation}\label{GJ.20b}
 \pi(\cM) \subset [m-C,m+C],\qquad \pi(\cN)\subset [n-C,n+C].
\end{equation}
\end{theorem}
Actually, $\cM$ consists of a single point, by a well-known argument using the convexity of $\Omega$. Also,
the uniqueness of the point $m$ follows by standard arguments from the convexity of $h$, while the uniqueness of $n$ is a general fact from Sturm theory.

A consequence of the theorem is that, while clearly the interval $I$ and the function $h$ do not determine $\Omega$ uniquely, these data do determine the location of the distinctive features of $u_1,u_2$ up to a bounded error, uniformly for all domains normalized as above (in particular, uniformly as $\abs{I}\to\infty$).

The question left open in \cite{Jer:DFNLCD,GriJer:SFECPD} is whether the result \eqref{GJ.20b} is
sharp in order of magnitude as $\abs{I}\to\infty$.  We will derive from Theorems \ref{GJ.4} and \ref{GJ.10} that this is in fact true:
\begin{theorem}\label{GJ.21a}
There is $c>0$ and, for each $N'>0$, a pair of domains $\Omega,\Omegatilde$, normalized as above and with $\pi(\Omega)=\pi(\Omegatilde)$ of length $N'$ and with the same height function, such that
\begin{equation}\label{GJ.22a}
\dist(\pi(\cM),\pi(\tilde{\cM}))> c,\qquad \dist (\pi(\cN),\pi(\tilde{\cN})) > c.
\end{equation}
\end{theorem}
See Figure \ref{figure2} at the end of the paper.
This should be put in contrast with the
main result of \cite{GriJer:AFNLCD} which states that the {\em length} of the interval $\pi(\cN)$ is not bounded away from zero, but actually
bounded above by $C/\abs{I}$, for an absolute constant $C$.

Our approach to Theorem \ref{GJ.4} is via the method of matched asymptotic expansions. This is carried out in Section \ref{asexp}. Here one has to deal with two limiting problems: The first one is the equation $(-\frac{d^2}{d\xi^2} + \lambda)\psi(\xi)=0$ on the interval $[0,1]$, and is easily solved explicitly. The second one is the equation $(\Delta+\pi^2)U=0$ on the unbounded domain $\Omega_\infty$; we need some results about generalized eigenfunctions on $\Omega_\infty$, which  are known from scattering theory. For completeness we give a direct
derivation of what we need in Section \ref{secinf}. Here, the quantity $a(\phi)$ arises. Theorem \ref{GJ.10} is proved in Section \ref{perturb} and Theorem \ref{GJ.20a} in Section \ref{apply}.

The asymptotic behavior of spectral quantities on degenerating spaces similar to the family $(\Omega_N)_{N>0}$ has been studied by many authors in different contexts. Regarding invariants involving all the eigenvalues (like the determinant of the Laplacian) we only mention
\cite{HasMazMel:ASAE}, \cite{Mue:EIMWB}
and \cite{ParWoj:STADZDDL}.
Our results are also related to the investigations of so-called thick graphs (or graph-like thin manifolds): When rescaling $\Omega_N$ by a factor $1/N$, one obtains a domain which is a $1/2N$-neighborhood of a unit interval, except for a fixed (scaled) perturbation at the ends of the interval. Instead of the interval (considered as a graph with two nodes and one edge of length one connecting them) one may consider more general embedded graphs, and their $\epsilon$-neighborhoods (or more general '$\epsilon$-thin' manifolds modeled on the graph). The convergence of eigenvalues, as $\epsilon$ tends to zero, to spectral data on the graph itself (then sometimes called a 'quantum graph') was studied in \cite{KucZen:CSMSCOG,RubSch:VPMCTSIBECLS,ExnPos:CSGLTM}, for the Neumann and closed problems. The Dirichlet and mixed boundary value problems are studied in the preprints \cite{Gri:SGNS} (where it is also proved that the asymptotic series constructed in this paper converge) and \cite{MolVai:LONTFSNT}, by different methods than the one used in this paper. The Dirichlet problem is more difficult to handle since the dependence on the counting parameter $k$ appears in a lower order term, cf. equation \eqref{GJ.5}. While in these papers the graph edges are always straight lines, the case of a curved line (but without the perturbation at the end) is considered in \cite{FreKre:LNSTCT}, where a nodal line theorem is proved.

\section{Eigenfunctions on the infinite domain} \label{secinf}
In this section we prove the results about generalized eigenfunctions  on $\Omega_\infty$ which will be
needed in the proof of Theorem \ref{GJ.4}. Since, for any $m$, the $m$th eigenvalue on $\Omega_N$ converges to $\pi^2$ as $N\to\infty$ (by domain comparison),
we need to consider the spectral value $\pi^2$, which is the bottom of the
continuous spectrum of $-\Delta_{\text{Dir}}$ on $\Omegainf$.

\begin{proposition}\label{GJ.1}
\mbox{}
\begin{enumerate}
\item
There is a unique function $U$ on $\Omega_\infty$ satisfying
\begin{align}
(\Delta+\pi^2)U &= 0 \quad \text{on } \Omega_\infty\notag\\
U &=0 \quad\text{at } \partial \Omega_\infty\label{GJ.1a}\\
U-x\sin \pi y & \ \text{ is bounded.}\notag
\end{align}

\item
For this function $U$ define
\begin{equation} \label{GJ.2}
a(\phi) = 2\int_0^1 U(0,y)\sin \pi y\, dy.
\end{equation}
Then
\begin{gather}\label{GJ.3}
U(x,y) = (x+a(\phi))\sin \pi y + r(x,y)\\
\intertext{where the remainder decays exponentially as $x\to\infty$, more precisely,}
 r(x,y) = \sum_{k=2}^\infty r^\park(x)\sin k\pi y,\quad
|d_x^p r^\park(x)| \leq C_p e^{-kx/2} \ \forall p,k.
\label{expdecay}
\end{gather}
Here, $C_p$ is independent of $k$ and $x$ and is bounded in terms of $\max \phi$.
\end{enumerate}
\end{proposition}

Note that in the special case where $\phi$ is a constant $A$,
the function $U$ is simply $(x+A)\sin\pi y$, so $a(\phi)=A$.
Therefore, for general $\phi$, the number $a(\phi)$ tells how much the 'standard' problem, with
$\Omega_\infty = [0,\infty)\times [0,1]$, has to be shifted so that
its first (generalized) eigenfunction coincides asymptotically with
that of the 'perturbed' problem. See also Remark \ref{remarkscattering}.

For simplicity, we assume all functions are real valued.
Basic to all considerations is the explicit form of solutions of the
homogeneous equation:

\begin{lemma}\label{GJ.90}
Assume $u$ solves $(\Delta+\pi^2)u=0$ in $x\geq x_0,\ y\in(0,1)$,
vanishes for $y=0$ and $y=1$ and has at most polynomial growth
as $x\to\infty$. Then
\begin{equation}\label{GJ.91}
u(x,y) = (A_1 + B_1x)\sin \pi y + \sum_{k=2}^\infty A_k e^{-\sqrt{k^2-1}\pi x}
\sin k\pi y  \quad\text{ for } x\geq x_0,
\end{equation}
for certain numbers $A_k$, $B_1$.
\end{lemma}
\begin{proof}
$u$ is smooth, so for each fixed $x$ it is the sum of its Fourier series,
$u(x,y) = \sum_{k=1}^\infty u_k(x) \sin k\pi y$, where
$u_k(x) = 2\int_0^1 u(x,y) \sin k\pi y\, dy$.
From $(\Delta+\pi^2)u=0$ one gets $u_k'' +(1-k^2)\pi^2 u_k = 0$. This gives
$u_1(x)=A_1+B_1x$ and $u_k(x)=A_ke^{-\sqrt{k^2-1}\pi x} +
B_ke^{\sqrt{k^2-1}\pi x}$ for $k\geq 2$.
Since $u$ is polynomially bounded, so is $u_k$, and therefore $B_k=0$
for $k\geq 2$.
\end{proof}

\begin{lemma}\label{GJ.100}
Let $w\in L^2(\Omega_\infty)$ be supported in $\Omega_1$.
\begin{enumerate}
\item
If $u\in H^1(\Oinf)$ solves
\begin{align}
(\Delta+\pi^2)u &= w \quad \text{on } \Oinf,
\notag\\
u &= 0\quad\text{at }\partial\Oinf,\label{GJ.100a}\\
u&\ \text{ is bounded }\notag
\end{align}
then
\begin{equation}\label{GJ.101}
\|u\|^2_{H^1(\Omega_1)} \leq -C (u,w)
\end{equation}
where $C$ is bounded in terms of the maximum of $\phi$.
\item
Problem
\eqref{GJ.100a} has  a unique solution $u$.
\end{enumerate}
\end{lemma}

\begin{proof}
(1)
We integrate by parts and use the support assumption to obtain
\begin{equation}\label{GJ.102}
 -(u,w) = -\int_{\Omega_1} u(\Delta+\pi^2)u
       = -\int_{\partial\Omega_1} uu_n
                      + \int_{\Omega_1}|\grad u|^2 - \pi^2\int_{\Omega_1} u^2
\end{equation}
where $u_n$ is the outward normal derivative.
Therefore, it is sufficient to prove the following two facts:
\begin{gather}
\text{If } u \text{ is as in \eqref{GJ.91} and bounded then }
\int_0^1 u(x_0,y)\frac{\partial u(x_0,y)}{\partial x}\,dy \leq 0.
\label{GJ.104}\\
\text{If } u\in H^1(\Omega_1)\text{ and } u=0 \text{ on }
\partial\Omega_1\setminus\{x=1\}\text{ then }
\int_{\Omega_1} |\grad u|^2 - \pi^2\int_{\Omega_1} u^2
\geq c \|u\|^2_{H^1(\Omega_1)}. \label{GJ.105}
\end{gather}
To prove \eqref{GJ.104}, observe that
$B_1=0$ in \eqref{GJ.91} since $u$ is bounded. Therefore,
$\int_0^1 u(x_0,y)\partial_x u(x_0,y) \,dy
= -\pi/2\sum_{k=2}^\infty A_k^2\sqrt{k^2-1} \exp(-2\sqrt{k^2-1}\pi x_0)\leq0$.

To prove \eqref{GJ.105}, consider the domain $\Otilde$ which is the union of
$\Omega_1$, the boundary piece $\{1\}\times(0,1)$, and the
reflection of $\Omega_1$ across this boundary.
Since $\Otilde\subset (-A,A)\times(0,1)$ for some
$A$, the first Dirichlet eigenvalue of $\Otilde$ is strictly bigger than
$\pi^2$, so it equals $(\pi^2+c)/(1-c)$ for some $c>0$.
The function on $\Otilde$ which equals $u$ for $x\leq 1$
and is symmetric with respect to the line $\{x=1\}$
is in $H^1_0(\Otilde)$, so we can
use it as test function and obtain
$(1-c)\int_{\Omega_1} |\grad u|^2 \geq (\pi^2+c)
\int_{\Omega_1} u^2$, which implies \eqref{GJ.105}. In this proof $A$ and therefore $c$ only depend on $\max\phi$.

(2)
Uniqueness is clear from 1.
To prove existence,
we reduce to a compact problem.
Define 'Dirichlet-to-Neumann operators' $N_1, N_\infty$, acting on functions
on $S=\{1\}\times (0,1)$, as follows:
Given $f\in H^1_0(S)\cap H^2(S)$, let $v_1,v_\infty$ be the solutions
of $(\Delta+\pi^2)v=0$ on $\Omega_1$, $D=(1,\infty)\times(0,1)$ respectively,
with boundary
values $f$ at $S$ and zero elsewhere, and $v_\infty$ bounded.
Existence and uniqueness follows for $v_1$ from the fact
that the first Dirichlet
eigenvalue of $\Omega_1$ is bigger than $\pi^2$, and for $v_\infty$
by explicit computation as in Lemma \ref{GJ.90}.
Set  $N_1 f := \partial v_1/\partial x_{|S}$ and
$N_\infty f := \partial v_\infty/\partial x_{|S}$. The restrictions exist
and are in $H^1(S)$ since $v_1,v_\infty$ are in $H^{5/2}$ near $S$
by standard regularity theory.

From \eqref{GJ.104} applied to $u=v_\infty$ we have
$(N_\infty f,f)\leq 0$, and from \eqref{GJ.105}
$(N_1f,f)=\int_{\partial\Omega_1} (v_1)_n v_1 =
\int_{\Omega_1} |\grad v_1|^2
-\pi^2 \int_{\Omega_1} v_1^2 \geq c\|v_1\|^2_{H^1(\Omega_1)} \geq \|f\|^2_{H^{1/2}(S)}$, so we obtain
\begin{equation}\label{GJ.103}
((N_1-N_\infty)f,f)\geq c\|f\|^2_{H^{1/2}(S)},\quad f\in H^1_0(S)\cap H^2(S).
\end{equation}
Along the same lines one sees that $(N_1f, g)\leq
\|f\|_{H^{1/2}(S)}\|g\|_{H^{1/2}(S)}$ for $f,g\in H^{1/2}_0(S)$,
and this shows that $N_1$ can be extended to a bounded operator
$H_0^{1/2}(S)\to H^{-1/2}(S)$, and similarly for $N_\infty$. By approximation,
\eqref{GJ.103} continues to hold for $f\in H^{1/2}_0(S)$.
This shows that $N_1-N_\infty$ has closed range and therefore is
surjective, for if $g\in H^{1/2}_0(S)$ is orthogonal to the range
then applying \eqref{GJ.103} to $g$  implies $g=0$.

Now we find a solution $u$ of \eqref{GJ.100a} as follows:
 Let $v$ be the unique solution of $(\Delta+\pi^2)v=w$
on $\Omega_1$, $v_{|\partial\Omega_1}=0$, and define $f$ on $S$ by
$(N_1-N_\infty)f = -\partial v/\partial x_{|S}$.
For $f$ find functions $v_1,v_\infty$ as above.
Define the function $u$ by $v+v_1$ on $\Omega_1$ and by
$v_\infty$ on $D$; at $S$, $u$ has the value $f$ from the left and right,
and $\partial u/\partial x$ is $\partial v/\partial x + N_1 f$
from the left and $N_\infty f$ from the right, which are equal by construction.
Therefore, $u$ is the desired solution.
\end{proof}

\begin{proof}[Proof of Proposition \ref{GJ.1}]
(1)
Let $\psi(x)$ be a smooth function vanishing for $x\leq 0$ and
equal to one for $x\geq 1$.
Let  $u$ be the solution of \eqref{GJ.100a} with
$w=-(\Delta+\pi^2)(\psi(x)x\sin \pi y)$. Then
$U=u+\psi(x)x\sin\pi y$ solves \eqref{GJ.1a}. From \eqref{GJ.101} we have $\norm{u}_{H^1(\Omega_1)}\leq C\norm{w}_{L^2(\Omega_1)}\leq C'$ and therefore also
\begin{equation}\label{GJ.103a}
\norm{U}_{H^1(\Omega_1)} \leq C,
\end{equation}
with $C$ only depending on $\max \phi$.
Uniqueness is clear from
the uniqueness for $u$.

(2) Apply Lemma \ref{GJ.90} to $U$ with $x_0=0$.
Since $U-x\sin\pi y$ is bounded, we have $B_1=1$.  Since, from \eqref{GJ.103a} and the trace theorem,
 $\sum_{k=1}^\infty A_k^2= \frac12\int_0^1 U(0,y)^2\,dy \leq C\norm{U}_{H^1(\Omega_1)}^2\leq C'$, this gives \eqref{GJ.3} and \eqref{expdecay}
with $a(\phi)=A_1$,
and evaluating $\int_0^1 U(0,y)\sin\pi y\,dy$
yields \eqref{GJ.2}.
\end{proof}

\begin{remark}
\label{remarkscattering}
We explain the relation of $a(\phi)$ to the scattering phase. Standard scattering theory (see, for example, \cite{Gui:TSQVB}) yields that  for $s>0$ close to zero  the equation $(\Delta + \pi^2 + s^2)E=0$, $E_{|\partial\Omega_\infty}=0$ has a unique polynomially bounded solution on $\Omega_\infty$ of the form
$$ E_s = (e^{-isx} + S(s)e^{isx})\,\sin \pi y + r_s(x,y)$$
for some number $S(s)$ (the {\em scattering matrix}) and a remainder $r_s(x,y)$ of the form \eqref{expdecay}. The function $S$ extends holomorphically to a neighborhood of zero in $\bbC$ and is real and of modulus one for real argument, hence may be written $S(s)=e^{i\gamma(s)}$ for a holomorphic function $\gamma$, the {\em scattering phase}. $r_s$ is also holomorphic in $s$ and the estimates \eqref{expdecay} are uniform in $s$ near zero, so one can take the limit $s\to 0$ to get a solution of $(\Delta+\pi^2)U=0$. This solution is bounded, hence constant equal to zero by Proposition \ref{GJ.1}. Therefore $S(0)=-1$ and $r_0\equiv 0$. One can get a nontrivial solution by taking $\lim_{s\to 0}\frac1s E_s = \frac d{ds}_{|s=0}E_s$, and this has leading term $-2i(x+\gamma'(0)/2)\sin\pi y$. Comparison with \eqref{GJ.3} then yields
\begin{equation}
\label{eqnscattphase}
a(\phi) = \frac12\gamma'(0).
\end{equation}
\end{remark}
\medskip

We will also need an extension of Proposition \ref{GJ.1}:
\begin{lemma}\label{scattlemma}
Assume $v$ is a smooth function on $\Oinf$
which vanishes at $\partial\Oinf$
and for $x\geq0$ has the form
$$
 v(x,y) = p(x)\sin\pi y + r(x,y),
$$
with $p$ a polynomial and $r$ satisfying the estimates \eqref{expdecay}.

Then any polynomially bounded solution of the problem
\begin{align}
\begin{split}
(\Delta+\pi^2)u &= v \quad \text{on } \Omega_\infty\\
u &=0 \quad\text{at } \partial \Omega_\infty
\end{split}
\label{problem}
\end{align}
has the same form $u(x,y)=q(x)\sin\pi y + s(x,y)$, where $s$ satisfies the estimates
\eqref{expdecay} and
\begin{equation}
q'' = p.
\label{qp}
\end{equation}
Such solutions exist, and are unique up to adding multiples of
$U$, where $U$ is defined in Proposition \ref{GJ.1}.
\end{lemma}
\begin{proof}
Let $u$ be a solution of \eqref{problem}. Taking the Fourier decomposition
of $u(x,\cdot)$ we get
$u(x,y) = q(x)\sin\pi y + \sum_{k=2}^\infty s^\park(x) \sin k\pi y$,
and then \eqref{problem} gives \eqref{qp} and, for each $k\geq2$,
\begin{equation}
(d_x^2 + 1-k^2)s^\park=r^\park.
\label{sequation}
\end{equation}
For any initial condition
$s^\park(0)=a_k$, \eqref{sequation} has the unique polynomially
bounded solution
\begin{equation}
s^\park(x) = a_k e^{-\alpha x}-\frac1{2\alpha}\int_0^\infty
\left(e^{-\alpha|x-z|}-e^{-\alpha(x+z)}\right) r^\park(z)\,dz
\label{specialsol}
\end{equation}
where $\alpha=\sqrt{k^2-1}$.
Since $u$ is given and $s^{(k)}(x)=2\int_0^1 u(x,y)\sin k\pi y\,dy$, we have
$a_k=2\int_0^1 u(0,y)\sin k\pi y\,dy$. An easy calculation shows
$|d_x^p s^\park(x)|\leq Ce^{-kx/2}$, and this proves the
first claim.

To prove existence of a solution, choose $q_0$ satisfying $q_0''=p$, and
a cutoff function $\chi$, equal to zero for $x\leq 1$ and
to one for $x\geq 2$, and set $u_0 (x,y) =\chi(x)\left( q_0(x)\sin\pi y +
\sum_{k=2}^\infty s^\park_0(x)\sin k\pi y\right)$. This is polynomially
bounded.
Then
$u=u_0+h$ solves \eqref{problem}
iff $h$ solves $(\Delta + \pi^2)h = v_1$, $h_{|\partial\Oinf}=0$,
where $v_1=v-(\Delta+\pi^2)u_0$ is compactly supported.  Lemma \ref{GJ.100}
guarantees the existence
of a bounded solution $h$.

Uniqueness is clear from Lemma \ref{GJ.90} and Proposition \ref{GJ.1}.
\end{proof}

\section{Asymptotic expansions of eigenfunctions and eigenvalues} \label{asexp}
In order to prove Theorem \ref{GJ.4} we use the idea of matched asymptotic expansions.
The strategy is this: First, we make reasonable guesses about the asymptotic behavior,
as $N\to\infty$, of the $m$th eigenvalue and of certain scaled limits of the eigenfunction. This leads to an ansatz in the form of formal asyptotic series in terms of powers of $1/N$, whose coefficients are undetermined numbers (for the eigenvalue) resp.\ functions (for the eigenfunction). The eigenvalue equation, the boundary condition and the condition that the various scaled limits must fit together ('match') in the transition region between different scaling regimes, yield a recursive system of equations for these coefficients. This system has a unique solution (Proposition \ref{ansatzworks}).  Given any approximation order, one then obtains a candidate for an approximate eigenvalue (by truncating the formal series), and also for an approximate eigenfunction, which is obtained by a suitable patching of the data from the different scaling regimes.
These candidates satisfy the eigenvalue equation with a small error (Proposition \ref{almostsol}), and from this we derive that they are close to actual eigenvalues. A domain comparison yields an a priori estimate on the actual eigenvalues, and this allows to conclude that all actual eigenvalues are obtained in this way, as well as a lower bound on the spectral gap. The spectral gap estimate then implies that the approximate eigenfunctions are close to actual eigenfunctions (Theorem \ref{ansatzcorrect}). Using this explicit information, it is easy to derive Theorem \ref{GJ.4}.

\subsection{The ansatz, formal eigenvalue and eigenfunction}
Fix an integer $m\geq 1$.
We want to find the $m$th eigenvalue and eigenfunction of $\Omega_N$,
asymptotically as $N\to\infty$. In this and the next subsection we simply write $\lambda$ for the $m$th eigenvalue  and $u$ for an associated
eigenfunction.
As a guide, recall that for the unperturbed case $\Omega_N = [0,N]\times
[0,1]$, we have (for $N>m/2$)
$$ \lambda = \pi^2 + N^{-2}m^2 \pi^2,\quad u(x,y) = \sin (N^{-1}m \pi x) \sin (\pi y).$$
Our ansatz is guided by the following expectations:
\begin{enumerate}
\item
The eigenvalue should have complete asymptotics:
\begin{align}
&\exists \lambda_{i}\in\RR,\quad i=2,3,4,\dots, \notag\\
\label{lexpansion}
&\lambda\sim \pi^2 + \sum_{i=2}^\infty N^{-i} \lambda_i,\quad N\to\infty.
\end{align}
Note that $\lambda=\pi^2+O(N^{-2})$ follows from domain comparison.
\item
At any fixed $x,y$, suitably normalized eigenfunctions should
converge as $N\to\infty$, and even have complete asymptotics:
\begin{align}
&\exists f_i : \Omega_\infty \to\RR,\quad i=0,1,2,\dots,
\notag\\
&u(x,y) \sim \sum_{i=0}^\infty N^{-i}f_i(x,y),\quad N\to\infty, \quad
(x,y)\in\Oinf.
\label{fexpansion}
\end{align}
\item
When fixing $\xi=x/N$ and $y$, and letting $N\to\infty$,
$u$ should converge, and even have complete asymptotics:
\begin{align}
&\exists g_j:[0,1]\times[0,1]\to\RR,\quad j=0,1,2,\dots,
\notag\\
&u(N\xi,y) \sim \sum_{j=0}^\infty N^{-j}g_j(\xi,y),\quad N\to\infty, \quad
(\xi,y)\in [0,1]\times[0,1].
\label{gexpansion}
\end{align}
\end{enumerate}
We get conditions on all the coefficients from three sources:
\begin{enumerate}
\item[(I)]
The equation
$(\Delta + \lambda) u = 0.$
Formally inserting the asymptotics above, differentiating term by term, and successively equating powers of
$N$, we get, for $i,j=0,1,2,\dots$:
\begin{align}
(\Delta + \pi^2) f_i &= - \sum_{k=2}^i \lambda_k f_{i-k}
\tag{f}\\
(\partial_y^2 + \pi^2) g_j &= -\partial_\xi^2 g_{j-2} - \sum_{k=2}^j \lambda_k
g_{j-k}
\tag{g}
\end{align}
(terms with negative indices are set equal to zero).
\item[(II)]
Boundary conditions on $\partial\Omega_N$.
\eqref{fexpansion}, \eqref{gexpansion} give
\begin{align}
f_i=0 &\text{ at }\partial \Oinf
\tag{bd f}\\
g_j=0 &\text{ at } \{y=0\}\cup\{y=1\}\cup\{\xi=1\}.
\tag{bd g}
\end{align}

We will prove below that (f), (bd f) imply that each $f_i$ has the form
\begin{equation}
f_i(x,y) = \varphi_i(x)\sin \pi y + r_i(x,y),\quad
\varphi_i(x)=\sum_{l=0}^i w_{il}\frac{x^l}{l!},
\label{fiform}
\end{equation}
with $r_i$ satisfying condition \eqref{expdecay}.
\item[(III)]
Matching conditions. To ensure small errors when patching the $f_i$ and the
$g_j$ to get an approximate eigenfunction,
we need to correlate the large $x$ behavior of the $f_i$ with the behavior
at $\xi=0$ of the $g_j$. This is done by formally writing
$ \sum_i N^{-i} f_i (x,y) = \sum_j N^{-j} g_j(x/N,y),$
expanding $f_i$ according to \eqref{fiform} and $g_j$ in Taylor series at
$\xi=0$ and equating the coefficients of $N^{-i}x^l$. This gives
$w_{il} \sin\pi y = \partial_\xi^l g_{i-l}(0,y)$.
This suggests to seek $g_j$ in the form
\begin{equation}
g_j (\xi,y) = \psi_j (\xi) \sin \pi y,
\label{gjform}
\end{equation}
and then the matching conditions read
\begin{equation}
w_{il} =  \frac{d^l \psi_{i-l}}{d\xi^l}(0),\quad l\leq i.
\tag{$M_{il} $}
\end{equation}

\end{enumerate}

Let us call a pair $(\sum_{i=0}^\infty N^{-i}f_i,\sum_{j=0}^\infty N^{-j}g_j)$
of formal series,
with $f_i$, $g_j$ of the form \eqref{fiform}, \eqref{gjform},
a {\em formal eigenfunction with formal eigenvalue}
$\pi^2+\sum_{i=2}^\infty N^{-i}\lambda_i$ if
(f), (g), (bd f),
(bd g) and ($M_{il}$) are satisfied for all indices, and not both
$f_0$, $g_0$ are identically
zero. Clearly, multiplying a formal eigenfunction
by a non-zero scalar, that is a series
$\sum_{i=0}^\infty N^{-i}a_i$ with $a_0\not=0$,
yields a formal eigenfunction again, with the same formal eigenvalue.

\begin{proposition}\label{ansatzworks}
If $\pi^2+\sum_{i=2}^\infty N^{-i}\lambda_i$ is a formal eigenvalue then
\begin{equation}
\lambda_2=m^2 \pi^2, \quad\text{ for some } m\in\NN.
\label{lambda2}
\end{equation}

Conversely, for each $m\in\NN$ there is a unique formal eigenvalue
with $\lambda_2=m^2 \pi^2$, and the formal eigenfunction is unique
up to multiplication by scalars.

Furthermore, we have
\begin{align}
\sum_{j=2}^\infty  N^{-j}\lambda_{j}&=\frac{m^2\pi^2}{(N+a)^2} + O(N^{-5}),\label{lambdaas}\\
\sum_{j=0}^\infty N^{-j}\psi_j(\frac{x}N) &= \sin m\pi\frac{x+a}{N+a} + O(N^{-3}),\label{psias}
\end{align}
where $a=a(\phi)$ is defined in \eqref{GJ.2}.

The $\lambda_j$, the coefficients of the $\varphi_j$, and the constants in the estimates \eqref{expdecay} of the remainders $r_j$ and in \eqref{lambdaas} and \eqref{psias} are all bounded in terms of $j$ and the maximum of $\phi$.
\end{proposition}

\begin{proof}
Because of \eqref{gjform} we may rewrite (g), (bd g) as
\begin{align}
(d_\xi^2 + \lambda_2)\psi_{j}& = -\sum_{l=3}^{j+2} \lambda_{l}\psi_{j+2-l},
\tag{$\psi$}\\
\psi_j(1)&=0,
\tag{bd $\psi$}
\end{align}
where we shifted the index by two.

To prove \eqref{lambda2},
note that $f_0$ is bounded and satisfies $(\Delta+\pi^2)f_0=0$
on $\Oinf$, hence is zero by
Lemma \ref{GJ.100}(2). Then, ($M_{00}$) gives $\psi_0(0)=0$, so we have
$$(d_\xi^2 + \lambda_2)\psi_0 = 0, \quad \psi_0(0)=\psi_0(1)=0,
\quad \psi_0 \not\equiv 0,$$
and this implies \eqref{lambda2}.

Now fix $m\in\NN$. We construct a formal eigenvalue and formal
eigenfunction with $\lambda_2=m^2 \pi^2$,
and satisfying the normalization condition
\begin{equation}
2\int_0^1 \psi_j(\xi)\sin m \pi \xi \, d\xi = \delta_{0j},\quad j=0,1,2,\dots,
\label{normalize}
\end{equation}
and simultaneously prove its uniqueness.
Since multiplying any formal eigenfunction by the scalar
$\left(2\sum_{j=0}^\infty N^{-j} \int_0^1 \psi_j(\xi)\sin m \pi \xi\,d\xi\right)^{-1}$
yields a formal eigenfunction satisfying \eqref{normalize}, this will prove
the Proposition.

First, by the argument proving \eqref{lambda2},
and by \eqref{normalize}, we must
have
$$ f_0 \equiv 0, \psi_0(\xi) = \sin m \pi \xi.$$
We now apply iteratively the following lemma.

\begin{lemma}
Let $J\geq 1$. Given $\psi_0,\dots,\psi_{J-1}$, $f_0,\dots,f_{J-1}$,
$\lambda_2,\dots,\lambda_{J+1}$ satisfying the equations
($\psi$), (f), ($M_{jl}$) and the boundary conditions for $j<J$,
there are unique
$\psi_J$, $f_J$, $\lambda_{J+2}$ satisfying these equations for $j=J$
and the normalization \eqref{normalize}.
\end{lemma}
\begin{proof}
First, we choose a solution $f_J$ of (f) (with $i=J$) of the form \eqref{fiform}, according to Lemma
\ref{scattlemma}, removing the indeterminacy by prescribing
$w_{J1} = \psi'_{J-1}(0)$.
This determines $w_{J0}$, and therefore $\psi_J(0)$
by $M_{J0}$. Next, equation ($\psi$) with $j=J$
has a solution with given values at 0 and 1, if and only
if the right hand side satisfies one linear condition, and then the solution
is unique up to multiples of $\psi_0=\sin m \pi y$, therefore uniquely determined by
condition \eqref{normalize}. The solvability
condition is obtained by taking the
scalar product of both sides with $\psi_0$ and integrating by parts on
the left. This gives
\begin{equation}
\psi_J(0) d_\xi\psi_0(0) = -\lambda_{J+2}/2,
\label{lambdaeqn}
\end{equation}
where we have used \eqref{normalize}, and this determines $\lambda_{J+2}$.

To finish the proof, we only need to check that $(M_{Jl})$ is satisfied
for $l\geq 2$. Now from \eqref{qp}, the polynomial $\varphi_J$
occuring in $f_J$
satisfies $\varphi_J'' = -\sum_{k=2}^J \lambda_k\varphi_{J-k}$. Equating
coefficients
of $\frac1{(l-2)!}x^{l-2}\sin \pi y$ we get $w_{Jl} = -\sum_{k=2}^J \lambda_k w_{J-k,l-2}$.
Using the matching conditions (with $j\leq J-2$) on the right and then the $(l-2)$th derivative
of equation ($\psi$), with $j=J-l$, we see that this sum equals
$d_\xi^l\psi_{J-l}(0)$.
\end{proof}
The boundedness of all quantities in terms of $\max\phi$ is also proved inductively, using the corresponding claims in Proposition \ref{GJ.1} and Lemma \ref{GJ.100}.

As an illustration, we carry this out for $J=1$: $\psi_0(\xi)=\sin m\pi\xi$ gives
$w_{11}=\psi_0'(0) = m\pi$, this determines $f_1=m\pi U$, and \eqref{GJ.3}
yields
$w_{10}=m\pi a$, hence $\psi_1(0)=m\pi a$. Equation \eqref{lambdaeqn}
now gives $\lambda_3=-2m^2\pi^2 a$, and we have
$\psi_1(\xi) = m\pi a (1-\xi)\cos m\pi a$.

It remains to check \eqref{lambdaas} and \eqref{psias}. The calculations of higher order terms
can be simplified by introducing new variables
$\xtilde=x+a$, $\Ntilde=N+a$, expressing all functions in terms
of $\xtilde$ and $\xitilde=\xtilde/\Ntilde$, and using formal series in
$\Ntilde$.
This must give the same result (after changing variables and up to
normalization) by uniqueness.
 We get in the $J=1$ step:
$\psitilde_0(\xitilde)=\sin m\pi \xitilde$ implies $\wtilde_{11}=m\pi$, so
$$\ftilde_1=m\pi\Utilde$$
(where $\Utilde(\xtilde,y)=U(x,y)$) as before, but now $\Utilde(\xtilde) = \xtilde \sin \pi y + O(e^{-\xtilde})$,
so $\wtilde_{10} = 0$ and thus $\psitilde_1(0)=0$, from which we get, using \eqref{lambdaeqn}
$$\psitilde_1\equiv 0,\lambdatilde_3=0.$$
The $J=2$ step yields $\wtilde_{21}=0$, and since $(\Delta+\pi^2)\ftilde_2=0$,
we get from Proposition \ref{GJ.1} that $\ftilde_2=0.$ As before, this gives
$\psitilde_2\equiv 0, \lambdatilde_4=0$, and this proves \eqref{lambdaas}, \eqref{psias}.
\end{proof}
\begin{remark}\label{remnextterm}
The next term can be obtained as follows.
For $J=3$ we get $\wtilde_{31}=0$, but now
$(\Delta+\pi^2)\ftilde_3 = -\lambdatilde_2\ftilde_1 = -m^3\pi^3 \Utilde$.
To obtain an expression for $\wtilde_{30}$, we multiply this with
$\Utilde$, integrate over $\Omegatilde_A$, apply Green's formula and let
$A\to\infty$. A short calculation (using \eqref{qp} also) gives
$$ \wtilde_{30} = 2m^3\pi^3 b,\quad b:= \lim_{A\to\infty} \left(
\int_{\Omegatilde_A}
\Utilde^2 d\xtilde dy - \frac{A^3}6\right),$$
and then \eqref{lambdaeqn} yields the first term missing in \eqref{lambdaas}
$$ \lambdatilde_5= -4m^4\pi^4 b.$$
\end{remark}

\subsection{Construction of an approximate solution from a formal solution}
For any order of approximation $M\in\NN$, we now use the formal solution obtained above to construct a candidate for an
approximate eigenvalue and eigenfunction on $\Omega_N$.
See Remark \ref{remmatching} for a motivation of our matching procedure.
In this subsection we still fix $m$ and omit it from the notation.

Choose cut-off functions $\chi_f$, $\chi_g$ on $\Omega_{N}$ as follows:
Choose a smooth function $\chi$ on $\RR$ which equals one on $(-\infty,1/2)$
and zero on $(3/4,\infty)$. Then set
$\chi_f(x,y) = \chi(x/N)$ and  $\chi_g (x,y) = 1-\chi(x)$.

Set
\begin{align}
f^\parM& = \sum_{i=0}^{M-1} N^{-i} f_i,\label{fM}\\
g^\parM(x,y)& = \sum_{j=0}^{M-1} N^{-j} g_j(\frac{x}N,y),\label{gM}\\
w^\parM(x,y)&= \sum_{i=0}^{M-1} N^{-i} \varphi_i(x) \sin \pi y.
\end{align}
$w^\parM$ describes the essential large $x$ behavior of $f^\parM$
and the small $\xi$
behavior of $g^\parM$.
Set
\begin{equation}
U^\parM = \chi_f f^\parM + \chi_g g^\parM - \chi_f \chi_g w^\parM.
\label{glue}
\end{equation}
That is, $U^\parM$ is given by $f^\parM$ for $x<1/2$, by $g^\parM$ for $x>(3/4) N$, and
by a smooth, appropriately scaled transition in between.
Finally, set
\begin{equation}
\label{lambdaM}
\lambda^\parM = \pi^2+\sum_{i=2}^{M-1} N^{-i} \lambda_i.
\end{equation}
\begin{proposition} \label{almostsol}
Denote $$V^\parM= (\Delta+\lambda^\parM)U^\parM.$$
For any $m,M,p\in\NN$ there are constants $c,C>0$ such that for all $N$
we have
\begin{align}
|\Delta^p V^\parM| &\leq C\left(\log N/N\right)^M,
\text{ uniformly in } \Omega_N,
\label{Uupperbound}\\
U^\parM, \Delta^p V^\parM &= 0 \text{ on }\partial \Omega_N,
\label{UVboundary}
\end{align}
and
\begin{equation}
\|U^\parM\|_{L^2(\Omega_N)} \geq c N^{1/2}.
\label{Ulowerbound}
\end{equation}
All constants $C$, as well as $c^{-1}$, are bounded in terms of $\max\phi$.
\end{proposition}

\begin{proof}
The idea is to split up \eqref{glue} in two ways: First, as
the $f$ term plus $g-w$, which gives estimates on the order of $(x/N)^M$
and then as the $g$ term plus $f-w$ which gives estimates
on the order of $N^{-M}+e^{-x}$. One of these is always bounded as in
\eqref{Uupperbound}.

Denote $\langle x\rangle=1+|x|$.
First, using the equations (f) for $i<M$ we get
$$ (\Delta+\lambda^\parM)f^\parM = \sum_{\substack{i,l\leq M-1\\ i+l\geq M}}
                 N^{-(i+l)} \lambda_l f_i.$$
Applying (f) again, we see that $\Delta f_i$ is a linear combination of
the $f_j$, $j\leq i$, and by induction over $p$ we get that
$\Delta^p(\Delta + \lambda^\parM) f^\parM$ is a linear combination of
$f_1,\dots,f_{M-1}$, with coefficients bounded by $N^{-M}$.
This implies
\begin{equation}
 \Delta^p(\Delta + \lambda^\parM)\chi_f f^\parM = O(\frac{\langle x\rangle^{M-1}}{N^M})
\label{festimate}
\end{equation}
uniformly in $\Omega_N$, since each $f_i=O(\langle x\rangle^i)$, and since derivatives
of $\chi_f$
are $O(N^{-1})$ and only occur where $x\geq N/2$, where each $N^{-i}f_i$, and
therefore $f^\parM$, and its
derivatives of any order are uniformly bounded.
Similarly, equation (g) and \eqref{gjform} give
$$ (\Delta+\lambda^\parM)g^\parM = \sum_{\substack{j,l\leq M-1\\ j+l\geq M+2}}
N^{-(j+l)} \lambda_j g_l,$$
and using (g) again and boundedness of the $g_j$ we get
\begin{equation}
\Delta^p (\Delta + \lambda^\parM)\chi_g g^\parM = O(N^{-M-2}) \quad
\text{ for }x\geq1.
\label{gestimate}
\end{equation}
Next, expanding $g_j(\xi,y)=\sum_{l=0}^{M-1} \partial_\xi^lg_j(0,y) \xi^l/l!
+ \xi^MR_j(\xi,y)$, with $R_j$ smooth, we get
\begin{multline*}
g^\parM - \chi_f w^\parM \\
= \sum_{j,l=0}^{M-1} N^{-j-l} x^l
\frac{\partial_\xi^lg_j(0,y)}{l!} - \chi(\frac{x}N)\sum_{i=0}^{M-1}\sum_{l\leq i}
N^{-i} x^l \frac{w_{il}}{l!} \sin\pi y + \xi^M \sum_{j=0}^{M-1} N^{-j}R_j.
\end{multline*}
Writing $i=j+l$ in the first sum, we see from ($M_{il}$) that all terms
of order at most $N^{-M+1}$ cancel, so we get
\begin{equation}
\Delta^p \chi_g(g^\parM - \chi_f w^\parM) = O(\frac{\langle x\rangle^M}{N^M}).
\label{gminusw}
\end{equation}
Finally, we have
\begin{equation}
\Delta^p \chi_f(f^\parM - \chi_g w^\parM) = O(e^{-x})\quad\text{for }x\geq 1
\label{fminusw}
\end{equation}
immediately from \eqref{fiform} and \eqref{expdecay}.

Writing $U^\parM = \chi_f f^\parM + \chi_g (g^\parM - \chi_f w^\parM)$, we
get from \eqref{festimate} and \eqref{gminusw}
\begin{equation}
\Delta^p(\Delta+\lambda^\parM)U^\parM = O(\frac{\langle x\rangle^M}{N^M}),
\label{Uest1}
\end{equation}
and writing $U^\parM = \chi_g g^\parM + \chi_f (f^\parM - \chi_g w^\parM)$
we get from \eqref{gestimate} and \eqref{fminusw}
\begin{equation}
\Delta^p(\Delta + \lambda^\parM)U^\parM = O (N^{-M-2} + e^{-x}),\quad x\geq 1.
\label{Uest2}
\end{equation}
Using \eqref{Uest1} for $x\leq M\log N$ and \eqref{Uest2} otherwise
we obtain \eqref{Uupperbound}.

The boundary conditions \eqref{UVboundary} are also clear from the arguments
above (note that $\chi_f,\chi_g$ only depend on $x$, so one only gets
$x$-derivatives of $f_i$, $g_j$ in the terms where the cut-offs are
differentiated).

\eqref{Ulowerbound} follows immediately from
the estimate
$|\sin m\pi x/N| \geq \sin\pi/8>0$, $x\in (N(1-\frac1{4m}),N(1-\frac1{8m}))$,
which implies $U^\parM \geq 1/10$ for these $x$ and $y\in (1/4,3/4)$, for
large $N$.
\end{proof}

\begin{remark} \label{remmatching}
Let us  clarify our procedure of obtaining asymptotic eigenfunctions,
by relating it to a simpler, 'compact' problem.
First recall how one may obtain a smooth
function $u(s,t)$ on $\RR^2$ with given Taylor expansions
$u\sim \sum_i t^i F_i(s)/i!$ at $t=0$ and $u\sim \sum_j s^j G_j(t)/j!$
at $s=0$, at least up to a certain order $M$:
First, such a $u$ exists iff the mixed derivatives of $u$ at
$(0,0)$ obtained from the two expansions agree, that is if
\begin{equation}
d_s^jF_i(0) = d_t^i G_j(0)\quad\text{ for all } i,j.
\label{matching}
\end{equation}
Calling this common value
$w_{ij}$ and setting, for a given order of approximation $M$,
$F^\parM =  \sum_{i=0}^{M-1} t^i F_i(s)/i!$,
$G^\parM = \sum_{j=0}^{M-1} s^j G_j(t)/j!$,
$w^\parM= \sum_{i,j=0}^{M-1} t^is^j w_{ij}/i!j!$,
one may set $$u^\parM = F^\parM + G^\parM - w^\parM.$$
From $$F^\parM-w^\parM = \sum_{i=0}^{M-1} \frac{t^i}{i!}
          \left(F_i(s) - \sum_{j=0}^{M-1} \frac{d_s^jF_i(0)}{j!} s^j \right)
          = O(s^M)$$
one sees that $u^\parM - G^\parM = O(s^M)$ uniformly for $(s,t)$ near zero,
and similarly $u^\parM - F^\parM = O(t^M)$ uniformly near zero, which was our
goal. These
estimates
continue to hold if one formally differentiates both sides any number of times.

This may be used to construct asymptotic solutions of partial differential equations: Let
$P$ be a partial differential operator of b-type, i.e.\ a polynomial
in $s\partial_s$, $t\partial_t$ with smooth coefficients. Suppose one can
determine the $F_i$ and $G_j$ so that $PF^\parM = O(t^M), PG^\parM =
O(s^M)$, uniformly near zero, (which amounts to solving a recursive set of ordinary
differential equations for the $F_i$ and the $G_j$) then
$$P u^\parM = P G^\parM + P(F^\parM - w^\parM)
= O(s^M)$$
and similarly $Pu^\parM = O(t^M)$, so
\begin{equation}
Pu^\parM = O(\min\{s^{M},t^{M}\}) = O((st)^{M/2}).
\label{opest}
\end{equation}

This relates to our problem as follows: We want to describe the eigenfunction
$u$ uniformly in $x$ and $N$, that is, as a function on
$$ D = \{(N,x,y):\, (x,y)\in\Omega_N\} \subset \RR^3.$$
In the sequel we suppress the $y$-dependence for simplicity.
Our ansatz postulates that $u$ has nice expansions in terms of smooth
functions of $x$ and $x/N$. This may be expressed as follows:
Introduce new variables $$s=\frac1x,\quad t=\frac{x}N$$ in the
subset $\{x\geq 1\}$ of $D$. Allowing the value $N=\infty$ (i.e.\ adding $\Omega_\infty$)
and then $s=0$, we get a
compactification $\Dtilde$  of $D$, given by adding a point at infinity for
each value of $t\in[0,1]$ and $y\in[0,1]$. What we prove is that $u$ extends to a function on
$\Dtilde$ which is smooth in $s$ and $t$ up to $s=0,t=0$.
The expansion \eqref{fexpansion} may be rewritten
$$ u \sim \sum_i \frac{t^i}{i!}F_i(s),\quad s\text{ fixed, where }
F_i(s) = s^if_i(\frac1s)i!,$$
and \eqref{gexpansion} becomes
$$ u\sim \sum_j \frac{s^j}{j!} G_j(t),\quad t\text{ fixed, where }
G_j(t) = t^jg_j(t)j!.$$
The matching conditions ($M_{il}$) are precisely the conditions
\eqref{matching}, with $l=i-j$.
Note that $F_i$ is smooth at $s=0$ by \eqref{fiform}.

The cutoff functions in \eqref{glue} must be introduced since
$F^\parM$ does not satisfy the boundary conditions at $t=1$
and $G^\parM$ does not satisfy the boundary conditions at $s\to\infty$,
i.e. at the left end of $\Omega_N$. In $(s,t)$ coordinates, the cut-offs
are simply functions of $s$ resp. $t$, and this motivates their choice
in \eqref{glue}.

We have $\partial_x = \frac{\partial s}{\partial x}\partial_s
+\frac{\partial t}{\partial x}\partial_t = -s^2\partial_s + st\partial_t,$
so the Laplacian is of b-type.
The estimate \eqref{Uupperbound} is actually stronger than what should
be expected from \eqref{opest} (which gives $O(N^{-M/2})$ only) since in our
problem the structure of $P$ yields $G_j=O(t^j)$ and
$F^\parM-w^\parM=O(e^{-1/s})$.
\end{remark}

\subsection{Closeness to actual solution, proof of Theorem \ref{GJ.4}}
\begin{theorem}\label{ansatzcorrect}
Denote the Dirichlet
eigenvalues of $-\Delta$ on $\Omega_N$ by $\mu_1<\mu_2\leq\dots$, and denote
the approximate $m$th eigenvalue constructed above by $\lambda^\parM_{m}$, and
the approximate eigenfunction by $U^\parM_{m}$.
Fix $m$. For sufficiently large
$N$ the first $m$ eigenvalues on $\Omega_N$ are simple, and for each $M$
\begin{equation}
|\mu_j - \lambda^\parM_{j}| = O(N^{-M}) \quad\text{for
  }j=1,\dots,m.
\label{evapprox}
\end{equation}
Furthermore, for $j\leq m$ there is
an eigenfunction $u_j$ for the eigenvalue $\mu_j$ satisfying, for any $\alpha\geq 0$,
\begin{equation}
 \sup_{\Omega_{N}} |D^\alpha (u_j-U^\parM_{j})| = O(N^{-M}),
\quad j=1,\dots,m.
\label{efapprox}
\end{equation}
The implied constants only depend on $M,j,\alpha$ and $\max\phi$.
\end{theorem}
\begin{proof}
Let $v_1,v_2,\dots$ be an orthonormal basis of eigenfunctions on $\Omega_N$,
corresponding to the eigenvalues $\mu_1,\mu_2,\dots$.
For fixed $j\in\{1,\dots,m\}$ write
\begin{equation}
U^\parM_{j} = \sum_l a_l v_l.
\label{Uexpansion}
\end{equation}
Then $a_l=(U_{j}^\parM,v_l)$ (scalar product in $L^2(\Omega_N)$).
For $V_{j}^\parM=(\Delta+\lambda_{j}^\parM)U_{j}^\parM$ we then
obtain $(V_{j}^\parM,v_l) =(\lambda_{j}^\parM-\mu_l)a_l$, using
integration by parts and $U^\parM_{j|\partial\Omega_N}=0$, and then by
induction $(\Delta^pV_{j}^\parM,v_l) = \mu_l^p(\lambda_{j}^\parM-\mu_l)a_l$ for
all $p\geq 0$ using \eqref{UVboundary}. Since \eqref{Uupperbound} implies
$\|\Delta^pV_{j}^\parM\|_{L^2}
\leq CN^{-M+1}$ we get
from Parseval's formula and \eqref{Ulowerbound}
\begin{align}
\sum_l |a_l|^2 &\geq cN
\label{ajlower}\\
\sum_l |a_l|^2 \mu_l^{2p}(\lambda_{j}^\parM-\mu_l)^2 &\leq C N^{-2M+2}.
\label{ajupper}
\end{align}

From \eqref{ajlower} and \eqref{ajupper}, with $p=0$, we get
$|\lambda_{j}^\parM-\mu_l|\leq CN^{-M+1/2}$ for some $l$.

Taking $M\geq4$, we get that
there is an eigenvalue $\mu_{l_j}$ in a $CN^{-3}$-neighborhood of
$\pi^2 + j^2\pi^2/N^2$, for each $j=1,\dots,m$. Since these neighborhoods are
disjoint for $N$ sufficiently large, we have
$\mu_{l_1}<\dots<\mu_{l_m}$, in particular $\mu_{l_m} \geq \mu_m$.
On the other hand, comparing $\Omega_N$ to the larger domain
$[-C,N]\times[0,1]$, we see that the $m$th eigenvalue (counting multiplicity)
of $\Omega_N$ is at
least $\pi^2+m^2\pi^2/(N+C)^2=\pi^2+m^2\pi^2/N^2 + O(N^{-3})$.
Therefore, $\mu_m\geq \mu_{l_m}$.
This implies that, for $N$ sufficiently large, $\mu_j=\mu_{l_j}$ for
$j=1\dots,m$ , so
the eigenvalues $\mu_j$ are all simple and
$|\mu_j - \lambda^\parM_{j}| = O(N^{-M+1/2})$.
Replacing $M$ by $M+1$ and subtracting $N^{-M}\lambda_{M,j}$
from $\lambda^{(M+1)}_j$ we obtain \eqref{evapprox}.

In particular, we have $|\mu_l-\lambda^\parM_{j}|\geq cN^{-2}$
for $l\not=j$, and
therefore
we get from \eqref{ajupper}
$ \sum_{l\not=j} \mu_l^{2p} |a_l|^2 \leq CN^{-2M+6}$
which means that
$$ \|\Delta^p(U^\parM_{j} - a_j v_j)\|_{L^2(\Omega_N)}
 \leq CN^{-M+3}.$$
By the Sobolev embedding theorem (applied to any unit width strip in $\Omega_{N}$)
we have $\|r\|_{C^\alpha(\Omega_{N}} \leq C\|\Delta^p r\|_{L^2(\Omega_N)} + C\|r\|_{L^2(\Omega_N)}$ for any
function $r$ on $\Omega_{N}$, whenever $2p>\alpha+1$.
Therefore, replacing
$M$ by $M+3$ and then subtracting terms of order $N^{-i}$, $i>M$, on the left
we get \eqref{efapprox} with $u_{j}=a_{j}v_{j}$.
\end{proof}

\begin{proof}[Proof of Theorem \ref{GJ.4}]
\eqref{GJ.5} follows with $M=5$ from \eqref{lambdaas}, \eqref{lambdaM} (where the index $m$ was omitted in the notation) and \eqref{evapprox} (for $j=m$).

For the eigenfunction we first recall \eqref{glue}, which gives for $x>1$ (where $\chi_g(x)=1$):
\begin{equation}\label{step1}
U_m^\parM = g^\parM + \chi_f(f^\parM - w^\parM).
\end{equation}
From \eqref{gM}, \eqref{gjform} and \eqref{psias} we have
\begin{equation}\label{step2}
 g^\parM(x,y) =  \sin m\pi\frac{x+a}{N+a}\,\sin \pi y + O(N^{-3}),
 \end{equation}
 and from \eqref{fminusw} we have for $x>3\log N$
 \begin{equation}\label{step3}
 \chi_f(f^\parM - w^\parM) = O(N^{-3}).
 \end{equation}
Clearly, the estimates \eqref{step2} and \eqref{step3} may be differentiated any number of times. Therefore,
 \eqref{efapprox}, \eqref{step1}, \eqref{step2} and \eqref{step3} give the eigenfunction
 estimate \eqref{GJ.6}.

Finally, from $[0,N]\times[0,1]\subset\Omega_N\subset [-\max\phi,N]\times[0,1]$ one has by domain monotonicity
$$ \pi^2 + \frac{m^2\pi^2}{(N+\max\phi)^2} \leq \mu_m(\Omega_N) \leq \pi^2 + \frac{m^2\pi^2}{N^2},$$
and combining this with \eqref{GJ.5} one obtains $0\leq a(\phi)\leq \max\phi$.
\end{proof}

\section{Perturbation of the domain} \label{perturb}
In this section we prove Theorem \ref{GJ.10}.

First, we derive an alternative formula for $a(\phi)$.
\eqref{GJ.2} can be rewritten
$a(\phi) = 2\int_{\partial\Omega_0} U\, \partial (x\sin\pi y)/\partial n$,
where $\partial/\partial n$ denotes differentiation in direction of the
outward unit normal $n$. Therefore,
by applying Green's formula on $\Omega_0$ and using
$(\Delta+\pi^2)U = (\Delta+\pi^2) (x\sin\pi y) = 0$
we obtain
\begin{equation}
\begin{split}
a(\phi) &= 2\int_{\partial\Omega_0}
           \frac{\partial U}{\partial n} x\sin\pi y \,ds(y) \\
&= 2\int_0^1 (\partial_x + \phi'(y)\partial_y)U_{|(-\phi(y),y)} \phi(y)
\sin\pi y\,dy
\label{aformula}
\end{split}
\end{equation}
since $n=(-1,-\phi')/\sqrt{1+(\phi')^2}$ and $ds=\sqrt{1+(\phi')^2}dy$.

Now fix $\phi$  and denote by $\Omega_N^\eps$ the domain
$\Omega_N$ defined using $\eps\phi$, and by $U^\eps$ the associated function
$U$ from Proposition \ref{GJ.1}.

Note that Theorem \ref{GJ.10} would follow from \eqref{aformula}
(with $\phi$ replaced by $\eps\phi$) if $U=U^\eps$ could be replaced by
$U^0(x,y) = x\sin\pi y$. Therefore, writing $v^\eps=U^\eps-U^0$ we only need to
show that
\begin{equation}
  \left\|\frac{\partial v^\eps}{\partial n}\right\|_{L^2(B^\eps)} =O(\eps),
\label{claim}
\end{equation}
where $B^\eps=\{(-\eps\phi(y),y):\,y\in[0,1]\}$ is the left boundary.
Since $v^\eps_{|B^\eps} = -\eps\phi \sin\pi y$, this follows from the
following lemma.

\begin{lemma}
Suppose $h$ is a function on $\partial\Omega_\infty$, supported in
$B=\{(-\phi(y),y):\, y\in[0,1]\}$, and $v$ solves
\begin{align*}
(\Delta+\pi^2)v &= 0 \quad \text{on } \Oinf,\\
v &= h\quad\text{at }\partial\Oinf,\\
v&\ \text{ is bounded }
\end{align*}
then
\begin{equation*}
\left\|\frac{\partial v}{\partial n}\right\|_{L^2(B)} \leq C\|h\|_{H^1(B)}
\end{equation*}
where $C$ is bounded in terms of the Lipschitz constant of $\phi$.
\end{lemma}
\begin{proof}
Write $v=u+H$, where $H$ is an extension of $h$ to $\Omega_\infty$, supported
in $\Omega_1$,  satisfying
$\|H\|_{H^1(\Omega_1)}\leq C\|h\|_{H^1(B)}$. Then $u$ satisfies
the assumptions of Lemma \ref{GJ.100} with $w=-(\Delta+\pi^2)H$,
so \eqref{GJ.101} gives $\|u\|_{H^1}\leq \|w\|_{H^{-1}}\leq
\|H\|_{H^1}$ and therefore
\begin{equation}
\|v\|_{H^1(\Omega_1)} \leq C \|h\|_{H^1(B)}.
\label{voinfest}
\end{equation}

Next, we choose a smooth cut-off function $\chi(x)$, equal to one in $x\leq 1/2$
and to zero in $x\geq 3/4$, and set $\vtilde = \chi v$. This satisfies
$\Delta \vtilde = w$, where $w:= -\pi^2\vtilde+2\grad\chi\grad v + (\Delta \chi)v$,
and $\vtilde_{|\partial\Omega_1} = h$, so standard estimates give

\begin{equation*}
\left\|\frac{\partial \vtilde}{\partial n}\right\|_{L^2(B)}
\leq C(\|w\|_{L^2(\Omega_1)} + \|h\|_{H^1(B)}) \leq C\|h\|_{H^1(B)}
\end{equation*}
using $\|w\|_{L^2}\leq C \|v\|_{H^1}$ and \eqref{voinfest}.
Since $v=\vtilde$ near $B$, this proves the lemma.
\end{proof}

\section{Maximum set and nodal line} \label{apply}

Here we prove Theorem \ref{GJ.21a}.
First, we obtain the following corollary of Theorem \ref{GJ.4}.
\begin{corollary}\label{coroll}
Consider the eigenfunctions $u_1$, $u_2$ on $\Omega_N$.
\begin{itemize}
\item[(a)] If $u_1$ assumes its maximum at a point $(x,y)$ then
\begin{equation}\label{maxloc}
|x - \frac{N-a(\phi)}2| = O(N^{-1}).
\end{equation}
\item[(b)]
If $(x,y)$ is an interior point of $\Omega_N$ with $u_2(x,y)=0$ then
\begin{equation}\label{nodloc}
|x-\frac{N-a(\phi)}2| = O(N^{-2}).
\end{equation}
\end{itemize}
\end{corollary}
\begin{proof} For shortness, we write $a=a(\phi)$.
(a) First, from \eqref{GJ.6} with $\alpha=0$ and from \eqref{GJ.6a} it follows that, at a maximum $(x,y)$, we must have $x\in[N/3,2N/3]$ and $\sin\pi y> 1/2$, for large $N$. Next, we use that $\partial_x u_1=0$ at a maximum.
From \eqref{GJ.6} with $\alpha=(1,0)$, i.e. taking $x$-derivatives,  we
obtain after multiplication by $N$ and division by $\sin\pi y$
$$\cos \pi\frac{x+a}{N+a} = O(N^{-2}).$$
With $x=\frac{N-a}2+\epsilon$ the expression on the left becomes $\cos (\frac\pi2 + \frac{\pi\epsilon}{N+a})
=-\sin\frac{\pi\epsilon}{N+a}$, so from $|\sin t|\geq |t|/2$ for small $t$ we get $|\epsilon| = O(N^{-1})$.

(b) First, by integrating \eqref{GJ.6} with $\alpha=(0,1)$ along the line from $(x,0)$ to $(x,y)$
we get $|u_m(x,y) - \sin m\pi\frac{x+a(\phi)}{N+a(\phi)}\,\sin\pi y| \leq C y N^{-3}$; by a similar
estimate in terms of distance to the upper boundary $y=1$ we obtain, after dividing by $\sin\pi y$, the improvement of \eqref{GJ.6},
\begin{equation}\label{GJ.6'}
 \sup_{x>3\log N} \left|\frac{u_m(x,y)}{\sin\pi y}-\sin m\pi\frac{x+a(\phi)}{N+a(\phi)}\right| = O(N^{-3}).
 \end{equation}
If $u_2(x,y)=0$ then, by Theorem 1 in \cite{Jer:DFNLCD}, we have $x\in[N/3,2N/3]$. Therefore, we obtain from \eqref{GJ.6'} (with $m=2$)
$$ \sin 2\pi\frac{x+a}{N+a} = O(N^{-3}).$$
As above this implies \eqref{nodloc}.
\end{proof}

We consider domains $\Omega_N$ of the form
\eqref{GJ.0} which are convex, i.e.\ with a concave function $\phi$.
By the corollary, to prove Theorem \ref{GJ.21a} it is enough to establish two concave functions $\phi$,
$\phitilde$ so that $a(\phi)\not=a(\phitilde)$ and the corresponding domains $\Omega_N$, $\Omegatilde_N$
have the same projection and height function.

\begin{figure}
\includegraphics{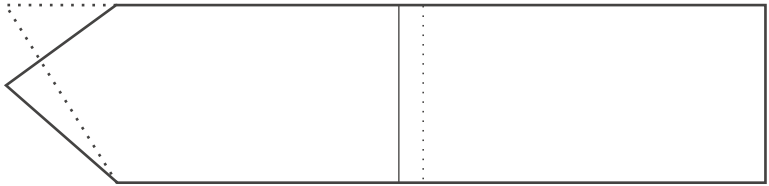}
\caption{Two domains with same height function but distinct location of the nodal line.}
\label{figure2}
\end{figure}

Let $\phi_0(y) = 1/2 - |y-1/2|$ and
$\phitilde_0(y) = y/2$.
See Figure \ref{figure2}.
Since $\phi_0$ is the symmetric decreasing
rearrangement of $\phitilde_0$ around the point $1/2$, and since $(\sin\pi
y)^2$ is symmetric decreasing itself, we have
$\int_0^1 \phi_0(y)(\sin\pi y)^2\,dy >
 \int_0^1 \phitilde_0(y)(\sin\pi y)^2\,dy$, so Theorem \ref{GJ.10}
implies $a(\epsilon\phi_0) \not= a(\epsilon\phitilde_0)$
for some sufficiently small
$\epsilon>0$. Since the domains associated with $\phi=\epsilon\phi_0$ and $\phitilde=\epsilon\phitilde_0$
clearly have the same height function, the theorem is proved.

\bibliographystyle{amsplain}
\bibliography{mypapers,dglib,mathlib}

\end{document}